\documentclass[12pt,a4paper]{article}
\usepackage{amsfonts,amsmath,amssymb,enumerate,amsthm,ifthen,graphpap,longtable,time,xspace,graphicx}
\usepackage[active]{srcltx}
\newcommand{\bbcmd}[1]{\mathop{\mathbb{#1}}\nolimits}
\newcommand{\rmcmd}[1]{\mathop{\mathrm{#1}}\nolimits}
\newcommand{\integer}{\bbcmd{Z}}
\newcommand{\real}{\bbcmd{R}}

\newcommand{\natur}{\bbcmd{N}}

\newcommand{\setcond}[2]{\left\{ #1 : #2 \right\}} 
\newcommand{\setcondbegin}[2]{\left\{ #1 : #2 \right.} 
\newcommand{\setcondend}[1]{\left.  #1 \right\}} 

\newcommand{\eps}{\varepsilon}

\newcommand{\card}[1]{| #1|}

\newcommand{\eop}{\hspace*{\fill}~$\square$} 
\theoremstyle{plain}
\newtheorem{theorem}{Theorem}[section]
\newtheorem*{theorem*}{Theorem}
\newtheorem{proposition}[theorem]{Proposition}
\newtheorem{lemma}[theorem]{Lemma}

\newtheorem*{corollary*}{Corollary}

\theoremstyle{definition}

\newtheorem{Remark}[theorem]{Remark}

\newcommand{\mycaption}[1]{\centering{\vspace{\medskipamount}\refstepcounter{figure}\textbf{Figure~\thefigure.} {#1}}}

\newenvironment{FigTab}[2]{
	\begin{figure}[htb]
	\setlength{\unitlength}{#2}
	\begin{center}
	\begin{tabular}{#1}
}{
    \end{tabular}
    \end{center}
    \end{figure}
}

\newcommand{\IncludeGraph}[2]{
	\includegraphics[#1]{{#2}}
}

\newcommand{\conv}{\rmcmd{conv}}

\title{Transversal Numbers over Subsets of Linear Spaces}
\author{\small G. Averkov\footnote{Institute of Mathematical Optimization, Faculty of Mathematics, University of Magdeburg, Universit\"atsplatz 2, 39106 Magdeburg, Germany, e-mail: averkov@math.uni-magdeburg} \ and R. Weismantel\footnote{Institute for Operations Research, ETH Z\"urich, R\"amistrasse 101, 8092 Z\"urich, Switzerland, e-mail: robert.weismantel@ifor.math.ethz.ch}}
\date{\small \today}
\begin{document} 
\maketitle

\begin{abstract} 
	Let $M$ be a subset of $\real^k$. It is an important question in the theory of linear inequalities to estimate the minimal number $h=h(M)$ such that every system of linear inequalities which is infeasible over $M$ has a subsystem of at most $h$ inequalities which is already infeasible over $M.$ This number $h(M)$ is said to be the Helly number of $M.$ In view of Helly's theorem, $h(\real^n)=n+1$ and, by the theorem due to Doignon, Bell and Scarf, $h(\integer^d)=2^d.$ We give a common extension of these equalities showing that $h( \real^n \times \integer^d) = (n+1) 2^d.$ We show that the fractional Helly number of the space $M \subseteq \real^d$ (with the convexity structure induced by $\real^d$) is at most $d+1$ as long as $h(M)$ is finite. Finally we give estimates for the Radon number of mixed integer spaces.
\end{abstract} 

\newtheoremstyle{itsemicolon}{}{}{\mdseries\rmfamily}{}{\itshape}{:}{ }{}
\newtheoremstyle{itdot}{}{}{\mdseries\rmfamily}{}{\itshape}{.}{ }{}
\theoremstyle{itdot}
\newtheorem*{msc*}{2000 Mathematics Subject Classification} 

\begin{msc*}
  Primary: 52A35,  90C11;  Secondary: 52C07
\end{msc*}

\newtheorem*{keywords*}{Key words and phrases}

\begin{keywords*}
Certificates of infeasibility; fractional Helly number; geometric transversal theory; Helly's theorem; integer lattice; linear inequalities; mixed integer programming; Radon's theorem
\end{keywords*}

\section{Introduction}

Let $M$ be a non-empty closed subset of $\real^k$. A subset $C$ of $M$ is said to be \emph{$M$-convex} if $C$ is intersection of $M$ and a convex set in $\real^k$. The set $M$ endowed with the collection of $M$-convex subsets becomes a convexity space (see \cite{MR1234493}). Let $h(M)$ be the \emph{Helly number} of $M$, i.e. the minimal possible $h$ satisfying the following condition.
\begin{itemize} 
	\item[$(H)$] Every finite collection $C_1,\ldots,C_m$ ($m \ge h$) of $M$-convex sets, for which every sub-collection of $h$ elements has non-empty intersection, necessarily satisfies $C_1 \cap \cdots \cap C_m \ne \emptyset.$ 
\end{itemize}
If $h$ as above does not exist we set $h(M) := \infty.$ The main purpose of this note is to study Helly type results in spaces $M \subseteq \real^k$ paying special attention to \emph{mixed integer spaces}, i.e. sets of the form $M=\real^n \times \integer^d$,  where $n, d \ge 0.$ For surveys of numerous Helly type results we refer to  \cite{MR0157289}, \cite{MR1242986},  \cite{MR1228043} and the monograph \cite{MR1234493}.

 It is known that $h(\real^n)=n+1,$ by the classical Helly Theorem (see \cite[Theorem\,1.1.6]{schneider}), and $h(\integer^d)=2^d$, by a result due to Doignon \cite[(4.2)]{MR0387090}, which was independently discovered by Bell \cite{Bell77} and  Scarf \cite{MR0452678} (see also  \cite[Theorem\,16.5]{schrijver} and \cite[p.\,176]{MR634767}). We obtain the following theorem. 

\begin{theorem} \label{mixed:helly:bounds} Let $M \subseteq \real^k$ be non-empty and closed and let $n,d \ge 0$ be integers. Then 
\begin{align} 
		 h( \real^n \times M) & \le (n+1) h(M), \label{mixed:with:reals} \\		
			2^d h(M) & \le h( M \times \integer^d ), \label{mixed:with:integers} \\
	h(\real^n \times \integer^d) & = (n+1) 2^d. \label{mixed:helly:eq}
\end{align}
\end{theorem}

Clearly, \eqref{mixed:helly:eq} is a direct consequence of  \eqref{mixed:with:reals}, \eqref{mixed:with:integers} and the theorems of Helly and Doignon. Equality \eqref{mixed:helly:eq} is a common extension of $h(\real^n)  = n+1$ and $h(\integer^d) = 2^d,$ and thus it can be viewed as Helly's theorem for mixed integer spaces. 

We notice that in general \eqref{mixed:with:integers} cannot be improved to equality (and thus, the case $M=\real^n$, for which we derive the equality, is quite likely an exception). For showing this we define $M:=\{0,1,2,2.5\}$ so that $h(M)=2.$ It turns out that $h(M \times \integer) \ge 5 > 2 h(M).$  In fact, consider the set $A:= \{(0,0),(1,0),(1,1),(2,1),(2.5,2)\}$, see also Fig.~\ref{example}. Then the five sets $A \setminus \{a\},$ $a \in A,$ are $(M \times \integer)$-convex and do not satisfy $(H)$ for $h=4,$ though each four of them have non-empty intersection.

\begin{FigTab}{cc}{0.7mm}
		\begin{picture}(70,65)
			\put(0,0){\IncludeGraph{width=70\unitlength}{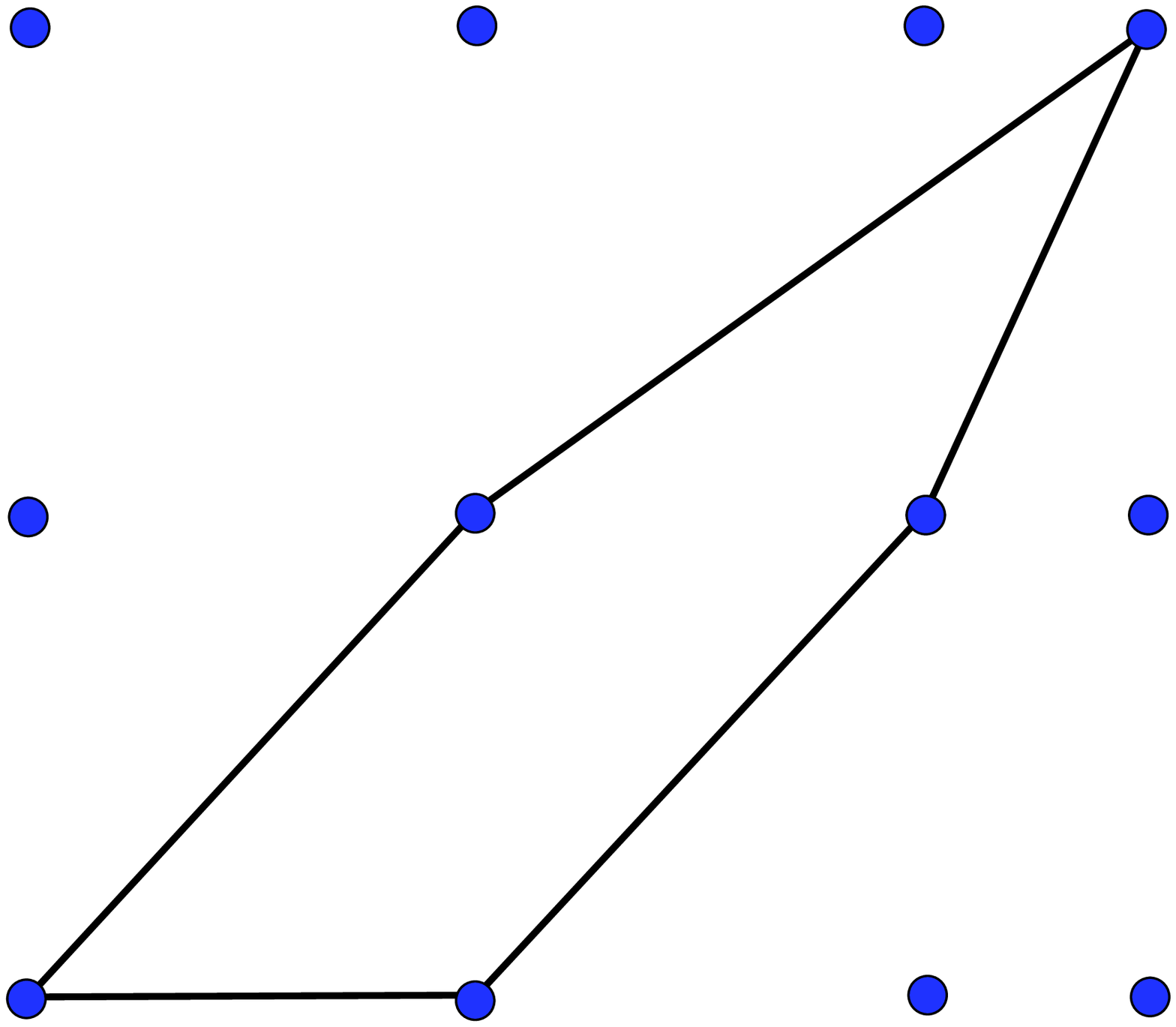}}
		\end{picture}
		\\
	\parbox[t]{0.85\textwidth}{\mycaption{Example for the case $h(M \times \integer^d) > 2^d h(M)$ for $d=1$\label{example}}}
\end{FigTab}

We wish to emphasize that theorems of Helly type are related, in a natural way, to the theory of linear inequalities. By this they also play a role in the theory of linear and integer programming, see for example \cite[Chapter\,7]{schrijver}, \cite[\S\,21]{MR1451876} and \cite{Clarkson95}. The following statement provides equivalent formulations of $h(M)$ in terms common for linear programming.

\begin{proposition} \label{helly:over:certificates}
	Let $M \subseteq \real^k$ be non-empty and $h \ge 0$ be an integer. Then (H) is equivalent to each of the following two conditions.
	\begin{itemize}
	\item[$(A)$] 	For every collection of affine-linear functions $a_1,\ldots, a_m$ $(m \ge h)$ on $\real^k$ either $a_1(x) \ge 0,\ldots, a_m(x) \ge 0$ has a solution $x \in M$ or, otherwise, there exist $1 \le i_1,\ldots, i_h \le m$ such that the system $a_{i_1}(x) \ge 0, \ldots, a_{i_h}(x) \ge 0$ has no solution $x \in M.$ 
	\item[$(B)$] For every collection of affine-linear functions $b_1,\ldots,b_m, c$ $(m \ge h-1)$ on $\real^k$ such that $\mu:=\sup \setcond{c(x)}{b_1(x) \ge 0,\ldots,b_m(x) \ge 0, \ x \in M}$ is finite there exist $1 \le i_1,\ldots,i_{h-1} \le m$ such that $\mu = \sup \setcond{c(x)}{b_{i_1}(x) \ge 0,\ldots,b_{i_{h-1}}(x) \ge 0, \ x \in M}.$ 
	\end{itemize}
\end{proposition}

Equivalence $(H) \Longleftrightarrow (B)$ above is an extension of the result of Scarf \cite{MR0452678} (see also \cite{Todd1977} and \cite[Corollary\,16.5a]{schrijver}).

Since $h(\real^n \times \integer^d)$ is linear in $n,$ formula \eqref{mixed:helly:eq} is in correspondence with the polynomial solvability of linear mixed integer programs in the case that the number of integer variables is fixed, see~\cite{lenstra83}. In terms of linear inequalities, \eqref{mixed:helly:eq} states that, in the mixed integer case, the (largest) number of inequalities in the insolvability certificate doubles if we introduce another integer variable and increases by $2^d$ (where $d$ is the number of integer variables) if we introduce another real variable.

\newcommand{\fh}{\bar{h}}
\newcommand{\aff}{\rmcmd{aff}}

We also wish to discuss fractional Helly's theorems in spaces $M \subseteq \real^k$. The \emph{fractional Helly number} $\fh(M)$ of $M$ is defined to be the minimal $h$ such that for every $0<\alpha \le 1$ there exists $0<\beta=\beta(\alpha,M) \le 1$ such that every collection of $M$-convex sets $C_1,\ldots,C_n$ ($n \in \natur$) satisfies the following condition. 

\begin{itemize} 
	\item[$(F)$] If $\bigcap_{i \in I} C_i \ne \emptyset$ for at least $\alpha \binom{n}{h}$ sets $I \subseteq \{1,\ldots,n\}$ of cardinality $h,$ then at least $\beta n$ elements of $C_1,\ldots,C_n$ have a common point. 
\end{itemize}

If $h$ as above does not exist, we set $\fh(M)=\infty.$  It is known that $\fh(\real^d)=\fh(\integer^d)=d+1$, where the relation for $\real^d$ is due to Katchalski and Liu  \cite{MR532152} and for $\integer^d$ due to B\'ar\'any and Matou\v{s}ek \cite{bar-mat-03}. B\'ar\'any and Matou\v{s}ek \cite[Remark at p.\,234]{bar-mat-03} point out that their arguments showing $\fh(\integer^d)=2^d$ do not use much of the geometry of $\integer^d$. Motivated by this observation, we prove the following extension.

\begin{theorem} \label{fh thm} 
	Let $M \subseteq \real^d$ be non-empty and closed and let $h(M)<\infty$. Then $\fh(M) \le d +1.$ 
\end{theorem}

Theorem~\ref{fh thm}  is a common extension of the fractional Helly results from \cite{MR532152} and \cite{bar-mat-03}. Clearly, Theorem~\ref{fh thm} applies to mixed integer spaces, and by this the fractional Helly number of $\integer^d \times \real^n$ is $d+n+1.$ The fractional Helly Theorem also implies the so-called $(p,q)$-theorem for $p \ge q \ge d+1$ in spaces $M \subseteq \real^d$ with finite Helly number $h(M)$, see \cite[p.\,229]{bar-mat-03} and \cite{MR1921545}.

The proof of Theorem~\ref{fh thm} follows the ideas from \cite{bar-mat-03}. In particular, we need a colored Helly's theorem for spaces $M.$ Let $I_1,\ldots,I_{d+1}$ be some pairwise disjoint sets of cardinality $t.$ We introduce the set $$K^{d+1}(t):= \setcond{ \{i_1,\ldots,i_{d+1}\}}{i_j \in I_j \ \mbox{for} \ j =1,\ldots,d+1},$$ the so-called \emph{complete $(d+1)$-uniform $(d+1)$-partite hyperpgraph}.  The foregoing notations are used in the following theorem.

\begin{theorem} \label{col h thm} {\upshape (Colored Helly's theorem for $M$-convex sets).}
	Let $M \subseteq \real^d$ be non-empty and let $h(M)<\infty.$ Then for every $r \ge 2$ there exists an integer $t=t(r,M)$ with the following property. For every family of $M$-convex sets $C_i$, $i \in I_1 \cup \ldots \cup I_{d+1}$ for which $\bigcap_{j=1}^{d+1} C_{i_j} \ne \emptyset$ for every $\{i_1,\ldots,i_{d+1}\} \in K^{d+1}(t),$ there exists $j \in \{1,\ldots,d+1\}$ and $R \subseteq I_j$ with $\card{R}=r$ such that $\bigcap_{i \in R} C_i \ne \emptyset.$ 
\end{theorem}

We do not need to give a proof of Theorem~\ref{col h thm}. In fact, the proof in our case is identical to the proof for the case $M=\integer^d$ from \cite{bar-mat-03}, based on Lov\'asz's colored Helly theorem, the theorem of Erd\H{o}s and Simonovits on super-saturated hypergraphs, and combinatorial arguments which rely only on the fact that the Helly number of $\integer^d$ is finite.

As a consequence of \eqref{mixed:helly:eq} we can also estimate the Radon number of mixed integer spaces. For a non-empty $M \subseteq \real^d$ the \emph{Radon number} $r(M)$ of the space $M$ is defined to be the minimal $r$ such that for every $A \subseteq M$ with $\card{A} \ge r$ there exist disjoint sets $B, C \subseteq A$ with $M \cap \conv B \cap \conv C  \ne \emptyset$. For $B$ and $C$ satisfying the previous condition the set $\{B,C\}$ is called a \emph{Radon partition} of $A$ in $M$. If $r$ as above does not exist we set $r(M) := \infty.$ By Radon's Theorem, $r(\real^n) = n+2$ (see \cite[Theorem\,1.1.5]{schneider}). It is known that $\Omega(2^d) = r(\integer^d) = O(d 2^d),$ see \cite{onn91}. The only known known value of $r(\integer^d)$ for $d \ge 2$ is $r(\integer^2) = 6$, see \cite{onn91} and \cite{MR2007959}. Furthermore, for $d=3$ one has $11 \le r(\integer^3) \le 17$  (see \cite{onn91} and \cite{MR2007959}). 

\begin{theorem} \label{radon:bounds}
	Let $M \subseteq \real^k$ be non-empty and closed and let $n, d\ge 0$ be integers. Then 
	\begin{align}
		&  r(M \times \integer^d) \ge (r(M) -1) 2^d  + 1, \label{radon:with:integers} \\
		 (n+1) 2^d  + 1 \le & \,  r(\real^n \times \integer^d) \le (n+d) (n+1) 2^d - n - d + 2. \label{radon:mixed:integers}
	\end{align}
\end{theorem}

We emphasize that the lower and upper bound in \eqref{radon:mixed:integers} are linear and quadratic in $n$, respectively. Thus, the exact asymptotics of $r(\real^n \times \integer^d)$ with respect to $n$ remains undetermined. Having \eqref{mixed:helly:eq} and \eqref{radon:mixed:integers}, it is suggestive to look for some type of Carath\'eodory's theorem in mixed integer spaces. However, the authors are currently not aware of any non-trivial notion of a Carath\'eodory number for $\real^n \times \integer^d.$ 

For dimension two $r(M)$ is uniquely determined by $h(M)$ with the only exceptional case $h(M)=4.$ 

\begin{theorem} \label{hel rad 2dim}
	Let $M \subseteq \real^2$ be non-empty. If $h(M) \ne 4,$ then $r(M) = h(M)+1.$ For $h(M)=4,$ one has $r(M) \in \{5,6\}$.
\end{theorem}

It is not hard to see that for $h(M)=4,$ both cases $r(M)=5$ and $r(M)=6$ are possible. It is known that $r(\integer^2)=6,$ and it is not hard to verify that $h(\integer \times \real)=5.$ It would also be interesting to study the relationship between $r(M)$ and $h(M)$ for $M \subseteq \real^d$ and $d \ge 3.$ 

\section{Proofs} 

In the proofs we use the standard terminology from the theory of polyhedra (see \cite{ziegler95}). A \emph{polyhedron} is the (possibly empty) intersection of finitely many closed  half-spaces. Bounded polyhedra are said to be \emph{polytopes}. If $P$ is a polyhedron, then faces of $P$ having dimension $\dim P -1$ are called \emph{facets}. By $\conv$ we denote the convex hull operation.

We first give the proof of Proposition~\ref{helly:over:certificates} since it is used as auxiliary statement in our main results.

\begin{proof}[Proof of Proposition~\ref{helly:over:certificates}]
	The implication $(H) \Longrightarrow (A)$ is trivial. Now, assuming that $(A)$ is fulfilled we derive $(H)$. Consider a collection $C_1,\ldots C_m$ ($m \ge h$) of $M$-convex sets such that every sub-collection of $h$ elements has non-empty intersection. For every $I \subseteq \{1,\ldots,m\}$ with $\card{I} = h$ we choose $p_I \in \bigcap_{i \in I} C_i.$ The polytope $$P_i := \conv \setcond{p_I}{I \subseteq \{1,\ldots,m\}, \ \card{I}=h, \ i \in I}$$ is a subset of $\conv C_i.$ Let $f_1,\ldots,f_{s}$ ($s \in \natur$) be affine-linear functions on $\real^n$ such that every $P_i$ can be given by $P_i = \setcond{x \in \real^n}{f_j(x) \ge 0 \ \text{for} \ j \in J}$ for an appropriate $J \subseteq \{1,\ldots,s\}$ and every $f_j$ is non-negative on some $P_i$. By construction, $\setcond{x \in M}{f_j(x) \ge 0 \ \text{for} \ j \in J} \ne \emptyset$ for every $J \subseteq \{1,\ldots,s\}$ with $\card{J}=h.$ Hence, in view of $(A)$,
	$$	
		\emptyset \ne \setcond{x \in M}{f_1(x) \ge 0,\ldots,f_s(x) \ge 0} \subseteq P_1 \cap \cdots \cap P_m \cap M \subseteq C_1 \cap \cdots \cap C_m. 
	$$

	Next, assuming that $(A)$ is fulfilled we derive $(B)$ by an argument analogous to the one given in \cite[Corollary\,16.5a]{schrijver}. Consider affine-linear functions $b_1,\ldots,b_m$ ($m \ge h-1$) such that the supremum $\mu$ in $(B)$ is finite. For every $t \in \natur$ the system $b_1(x) \ge 0,\ldots, b_m(x) \ge 0, c(x) \ge \mu+1/t$ has no solution $x \in M,$ and hence its subsystem consisting of $h$ inequalities inequalities has no solution in $M.$ Each such subsystem contains the inequality $c(x) \ge \mu+1/t.$ Thus, it follows that there exist $i_1,\ldots, i_{h-1}$ such that the system $b_{i_1}(x) \ge 0,\ldots, b_{i_{h-1}}(x) \ge 0, c(x) \ge \mu+1/t$ has no solution $x \in M$ for infinitely many $t.$ The latter obviously implies $(B).$ 

	In order to show $(B) \Longrightarrow (A)$ we assume that $(A)$ is not fulfilled and derive that $(B)$ is not fulfilled, as well. Let $a_1,\ldots,a_m$ $(m \ge h)$ be affine-linear functions such that $\setcond{x \in M}{a_1(x) \ge 0,\ldots, a_m(x) \ge 0} = \emptyset$ and for every $I \subseteq \{1,\ldots,m\}$ such that $\card{I}=h$ there exists a $p_I \in \setcond{x \in M}{a_i(x) \ge 0 \ \mbox{for} \ i \in I}.$ Without loss of generality we may assume that every subsystem of $a_1(x) \ge 0,\ldots, a_m(x) \ge 0$ consisting of $m-1$ inequalities is solvable over $M.$ In fact, otherwise we can redefine the system $a_1(x) \ge 0,\ldots,a_m(x) \ge 0$ by passing to a proper subsystem. Consider $c(x) := a_m(x)$ and $b_j(x) := a_j(x)$ for $j \in \{1,\ldots,m-1\}.$ Choose affine-linear functions $b_m(x),\ldots,b_{m+k}(x)$ such that $\setcond{x \in \real^k}{b_m(x) \ge 0,\ldots,b_{m+k}(x) \ge 0}$ is a simplex which contains all $p_I$'s introduced above.  Then 
	\[-\infty<\sup \setcond{c(x)}{b_1(x) \ge 0,\ldots,b_{m+k}(x) \ge 0, \ x \in M} < 0.\]
Furthermore, for all $1 \le i_1,\ldots,i_{h-1} \le m+k$ one has 
\[\sup \setcond{c(x) }{b_{i_1}(x) \ge 0,\ldots, b_{i_{h-1}}(x) \ge 0, x \in M} \ge c(p_I) \ge 0\] 
for every $I \subseteq \{1,\ldots,m\}$ satisfying $\card{I} = h$ and $\left(\{i_1,\ldots,i_{h-1}\} \cap \{1,\ldots,m-1\} \right) \cup \{ m \} \subseteq I.$ Hence $(B)$ is not fulfilled (for $m+k$ in place of $m$), and we are done.
\end{proof}

\begin{lemma} \label{helly:over:cert:lem}
	Let $M \subseteq \real^k$ be non-empty and closed and let $h \in \natur,$ $h \ge k+1.$ Then $(H)$ is equivalent to the following condition.
	\begin{itemize}
		\item[$(A')$]  For every choice of affine-linear functions $a_1,\ldots,a_m$ ($m \ge h$) on $\real^k$ such that the polyhedron 
	\begin{equation} \label{P:over:a}
		P:= \setcond{x \in \real^k}{a_1(x) \ge 0,\ldots,a_m(x) \ge 0}
	\end{equation} 
	 satisfies the conditions:
	\begin{enumerate} 
	\item \label{P:bounded} $P$ is bounded,
	\item \label{P:full:dim} $P$ is $k$-dimensional,
	\item \label{PdisjM} $P \cap M = \emptyset$,
	\end{enumerate}
	one necessarily has 
	\begin{equation} \label{certificate} 
		\setcond{x \in M}{a_{i_1}(x) \ge 0,\ldots,a_{i_h}(x) \ge 0} = \emptyset
	\end{equation}
	for some $1 \le i_1,\ldots,i_h \le m.$
	\end{itemize}
\end{lemma}
\begin{proof} 
	The implication $(H) \Longrightarrow (A')$ is trivial. Now, assume that $(A')$ is fulfilled.  We will show that $(A')$ also holds when we drop out Conditions~\ref{P:bounded} and \ref{P:full:dim} (which, in view of Proposition~\ref{helly:over:certificates}, yields the sufficiency). Assume that $P$ given by \eqref{P:over:a} satisfies Conditions~\ref{P:bounded} and \ref{PdisjM} but does not satisfy Condition~\ref{P:full:dim}. If $P= \emptyset,$ the existence of $i_1,\ldots,i_h$ satisfying  \eqref{certificate} follows from Helly's theorem for $\real^k.$ Assume that $P \ne \emptyset.$ Then, employing the closedness of $M$ and compactness of $P$, we see that there exists an $\eps>0$ such that $P_\eps := \setcond{x \in \real^k}{a_1(x) + \eps \ge 0,\ldots, a_m(x) + \eps \ge 0}$ is a $k$-dimensional polytope with $P_\eps \cap M = \emptyset.$ Consequently, applying $(A')$ for the affine-linear functions $a_1(x) +\eps, \ldots, a_m(x)+ \eps,$ we obtain $$\setcond{x \in \real^k}{a_{i_1}(x) + \eps \ge 0,\ldots,a_{i_h}(x) + \eps \ge 0} \cap M = \emptyset$$ for some $1 \le i_1,\ldots,i_h \le m.$ The latter implies \eqref{certificate}. Thus, $(A')$ still holds when we drop out Condition~\ref{P:full:dim}. Take $P$ given by \eqref{P:over:a} which satisfies Condition~\ref{PdisjM} but does not satisfy Condition~\ref{P:bounded}. For $t \in \natur$ we introduce the polytope 
	
	$$Q_t := \setcond{x \in \real^k}{x \in P, \ \pm x_1+t \ge 0,\ldots, \pm x_k + t \ge 0}$$ 
	where $x_1,\ldots,x_k$ are coordinates of $x$. Applying $(A')$ (with dropped out Condition~\ref{P:full:dim}) for the affine-linear functions $\pm x_i + t$ ($i \in \{1,\ldots,k\}$), $a_j(x)$ ($j \in \{1,\ldots,m\}$) that define $Q_t$ we find sets $I \subseteq \{1,\ldots,m\}$, $J^+, J^- \subseteq \{1,\ldots,k\}$ (a priori depending on $t$) such that $\card{I}+\card{J^+} + \card{J^-} \le h$ and 
	\begin{align} 
		\setcondbegin{x \in  M}{} & \, a_i(x) \ge 0 \ \mbox{for} \  i \in I, \nonumber \\
		& \setcondend{x_j+t \ge 0 \ \mbox{for} \ j \in J^+, \ -x_j+t \ge 0 \ \mbox{for} \ j \in J^-}  = \emptyset. \label{IJ:cond}
	\end{align}
	Since the index sets $\{1,\ldots,m\}, \{1,\ldots,k\}$ are finite we can fix $I, J^-, J^+$ independent of $t$ and such that \eqref{IJ:cond} holds for infinitely many $t$'s. Then $\setcond{x \in M}{a_i(x) \ge 0 \ \mbox{for} \ i \in I}  = \emptyset$, $\card{I} \le h,$ and the assertion follows.
\end{proof}

\begin{proof}[Proof of Theorem~\ref{mixed:helly:bounds}] Inequalities \eqref{mixed:with:reals} and \eqref{mixed:with:integers} are trivial if $h(M) = \infty.$ Thus, we assume $h(M)< \infty.$ Furthermore, without loss of generality we assume that $M$ affinely spans $\real^k$ so that $h(M) \ge k+1.$ We derive \eqref{mixed:with:reals} with the help of Lemma~\ref{helly:over:cert:lem}. Consider arbitrary affine-linear functions $b_1,\ldots,b_s$  ($s \in \natur$) on $\real^n \times \real^k$ such that $$P:= \setcond{(x,y) \in \real^n \times \real^k}{b_1(x,y) \ge 0,\ldots, b_s(x,y) \ge 0}$$ is $(n+k)$-dimensional, bounded and $P \cap (\real^k \times M) = \emptyset$. Let $T$ be the canonical projection from $\real^n \times \real^k$ onto $\real^k.$ The $k$-dimensional polytope $T(P)$ can be represented by $$T(P) = \setcond{y \in \real^k}{a_1(y),\ldots, a_{m}(y) \ge 0},$$ where $a_1,\ldots,a_{m}$ are affine-linear functions on $\real^k$ such that for each $j \in \{1,\ldots,m\},$ the set $F_j:= \setcond{y \in T(P)}{a_j(y) = 0}$ is a facet of $T(P).$ Hence $G_j:=T^{-1}(F_j) \cap P$ is a face of $P$ of dimension at least $k -1$. Consequently, the cone $N_j$ of affine-linear functions $f(x,y)$ vanishing on $G_j$ and non-negative on $P$ has dimension at most $(n+k)-(k-1) = n+1.$ The cone $N_j$ is generated by those $b_i$, $i \in \{1,\ldots,s\}$, which vanish on some facet of $P$ that contains $G_j.$ The function $a_j(y)$, $y \in \real^k$, can also be viewed as an affine-linear function on $\real^n \times \real^k.$ Thus, by Carath\'edory's Theorem for convex cones (cf. \cite[\S\,7.7]{schrijver}) applied to the function $a_j(y)$ in the cone $N_j$, there exists $I_j \subseteq \{1,\ldots,s\}$ such that $\card{I_j} \le n+1$ and 
\begin{equation} \label{a_j:over:b_i}
	a_j(y) = \sum_{i \in I_j} \lambda_{i,j} b_i(x,y)
\end{equation} 
 for appropriate $\lambda_{i,j} \ge 0$ ($i \in I_j$). By the definition of $h(M),$ there exists $J \subseteq \{1,\ldots,m\}$ with $\card{J} \le h(M)$ such that $\setcond{y \in M}{a_j(y) \ge 0 \ \mbox{for} \ j \in J}  = \emptyset.$ It follows that  the system $b_i(x,y) \ge 0$  with $i \in \bigcup_{j \in J} I_j$  has no solution $(x,y) \in \real^n \times M$. This system consists of at most $(n+1) h(M)$ inequalities. Hence, in view of Lemma~\ref{helly:over:cert:lem}, we arrive at \eqref{mixed:with:reals}

	Let us show \eqref{mixed:with:integers}. Let $h:=h(M)$ and $C_1,\ldots,C_h$ be $M$-convex sets such that $C_1 \cap \ldots \cap C_h = \emptyset$ but every sub-collection of $C_1,\ldots,C_h$ consisting of $h-1$ elements has non-empty intersection. Then the collection $(C_i \times \{j\}) \cup (C_i \times \{1-j\})$, where $i \in \{1,\ldots,h \}$, $j \in \{ 0,1\}$ consists of $2 h$ $(M \times \integer)$-convex sets, has empty intersection and the intersection over  every of its proper sub-collections is non-empty. This shows the bound $h(M \times \integer) \ge 2 h(M).$ The general bound $h(M \times \integer^d) \ge 2^d h(M)$ is obtained by induction on $d.$

	Equality \eqref{mixed:helly:eq} is a consequence of \eqref{mixed:with:reals}, \eqref{mixed:with:integers} and the theorems of Helly and Doignon.
\end{proof}

We remark that representation \eqref{a_j:over:b_i} can be viewed as Farkas type certificate of insolvability of a system of linear inequalities, see also \cite[\S\,13.1]{BertWeiBook} and \cite{CertificatesMixed2008} for related results.

Next we work towards the proof of Theorem~\ref{fh thm}. We show that the main tools in the proof of the fractional Helly theorem for $\integer^d$ given in \cite{bar-mat-03} can also be applied for sets $M \subseteq \real^d$ with a finite Helly number.  First we obtain a weak form of the fractional Helly theorem. 

\begin{theorem} \label{weak fh thm}
	Let $M \subseteq \real^d$ be non-empty and closed and let $h(M)<\infty.$ Then $\fh(M) \le h(M).$ 
\end{theorem}
\begin{proof} 
	\newcommand{\II}{\mathcal{I}} We slightly adjust the proof from \cite[Proof of Theorem~2.5]{bar-mat-03}.
	Let $h:=h(M).$ 
	Consider $M$-convex sets $C_1,\ldots,C_n$ ($n \in \natur$). Let $\II$ be the collection of those $h$-element subsets $I$ of $\{1,\ldots,n\}$ for which $\bigcap_{i \in I} C_i \ne \emptyset.$ Assume that $\card{\II} \ge \alpha \binom{n}{h}$ for some $0 \le \alpha < 1.$ For every $I \in \II$ we choose $p_I \in \bigcap_{i \in I} C_i$ and introduce the polytopes $P_i := \conv \setcond{ p_I}{I \in \II, \ i \in I},$ $i=1,\ldots,n.$ By construction, $P_i \subseteq C_i$ for every $i,$ and $p_I \in \bigcap_{i \in I} P_i.$ If $S \subseteq \{1,\ldots,n\}$ we shall write $P_S:= \bigcap_{i \in S} P_i.$  It is known that for a given non-empty, compact convex set $K$ and almost all directions $u$, the direction $u$ is an outward normal of precisely one boundary point of $K$; for a precise formulation see \cite[Theorem\,2.2.9]{schneider}. Applying this result to the sets $\conv (P_S \cap M)$ with $S \subseteq \{1,\ldots,n\}$, we see that there exists an affine function $a$ such that $a$ is maximized in exactly one point $x_S$ on $P_S \cap M$ (as long as $P_S \cap M$ is non-empty).
	
	\newcommand{\QQ}{\mathcal{Q}}

	For every $I \in \II$ there exists an $(h-1)$-element subset $J=J(I)$ of $I$ such that $x_{J}=x_{I}$ is the unique point maximizing $a(x)$ for $x \in M \cap P_J.$ In fact, for $H_{I} := \setcond{x \in M}{a(x) > a(x_I)}$ the family $\setcond{P_i \cap M}{i \in I} \cup \{H_{I}\}$ has empty intersection. Therefore, by the definition of $h(M),$ some $h$-element subfamily of this family has empty intersection. The elements of this subfamily which do not coincide with $H_{I}$ determine $J.$ 

	There are at most $\binom{n}{h-1}$ possible sets $J$ and at least $\alpha \binom{n}{h}$ different sets $I$. Thus, for a suitable $\beta=\beta(\alpha,h)>0$, some $J=:J^\ast$ is assigned to at least $\beta n$ different sets $I$. Each such $I$ has exactly one $i \not\in J^\ast$, and $x_{J^\ast}$ is a common point of these (at least $\beta n$ many) sets $P_i$.
\end{proof}

\begin{proof}[Proof of Theorem~\ref{fh thm}]
	The proof is a consequence of Theorems~\ref{col h thm}, \ref{weak fh thm} and the Erd\H{o}s-Simonovits theorem on super-saturated hypergraphs; for details see \cite[p.\,232]{bar-mat-03}. 
\end{proof}

\begin{proof}[Proof of Theorem~\ref{radon:bounds}]
	We exclude the trivial case $r(M) = \infty.$ The inequality $r(M \times \integer) \ge 2 r(M) -1$ can be shown following the idea from \cite[proof of Proposition~2.1]{onn91} (see also \cite[pp.\,176-177]{MR1234493}). Consider a set $A \subseteq M$ with $\card{A} = r(M)-1$ such that $M \cap \conv B \cap \conv C = \emptyset$ for all disjoint $B, C \subseteq A.$ Then $(M \times \integer) \cap \conv B \cap  \conv C = \emptyset$ for all disjoint $B, C \subseteq A \times \{0,1\}$. Thus $r(M \times \integer) \ge \card{A \times \{0,1\}} +1 = 2 r(M) -1.$ The general bound \eqref{radon:with:integers} is obtained by induction on $d.$ 

	The lower bound in \eqref{radon:mixed:integers} is a direct consequence of \eqref{radon:with:integers} and Radon's Theorem. The upper bound in \eqref{radon:mixed:integers} follows from the known inequality $r(M) \le k (h(M)-1) + 2$ (see \cite{MR0514026}, \cite[p.\,169]{MR1234493}) and \eqref{mixed:helly:eq}.
\end{proof}

\begin{proof}[Proof of Theorem~\ref{hel rad 2dim}]
	The case $h(M) = \infty$ is trivial. It is known that $r(M) \ge h(M)+1.$ For $h(M) \le 3$ it is easy to establish the inequality $r(M) \le h(M)+1.$ For the case $h(M)=4$ we show $r(M) \le 6$ following the idea from \cite[p.\,182]{MR2007959}. Assume the contrary, there exists a six-point set $A \subseteq M$ which does not possess a Radon partition in $M.$ Then $\conv A$ is a hexagon. We notice that any four of the six sets $A \setminus \{a\}, \ a \in A$ have non-empty intersection. Thus, by the definition of $h(M),$ $\emptyset \ne \bigcap_{a \in A} \conv (A \setminus \{a\}) \cap M = \conv \{a_1,a_3,a_5\} \cap \conv \{a_2, a_4, a_6\} \cap M,$ where $a_1,\ldots,a_6$ are consecutive vertices of $\conv A.$ Hence $\left\{ \{a_1,a_3,a_5\},\{a_2,a_4,a_6\} \right\}$ is a Radon partition of $A$ in the space $M,$ a contradiction. Now we consider the case that $h:=h(M) \ge 5$ and show that $r(M) \le h(M)+1.$ Let $A \subseteq M$ be a set of cardinality $h+1.$ If $A$ is not a vertex set of a convex polygon, it possesses a Radon partition. Assume that $A$ is a vertex set of a convex polygon. Then, by the definition of $h(M),$ we can choose 
	\begin{equation} \label{p in core} 
		p \in \bigcap_{a \in A} \conv (A \setminus \{a\}) \cap M. 
	\end{equation} 
	By Carath\'eodory's theorem (for $\real^2$) $p$ is in $\conv T$, for some three-element subset $T$ of $A.$ If $\conv T$ and $\conv A$ do not share edges, then taking into account \eqref{p in core} and the fact that $A$ is a vertex set of $\conv A$ we get $p \in \conv (A \setminus T) \cap \conv T \cap M.$ Consider the case that $\conv T$ and $\conv A$ share an edge. First notice that $\conv A$ and $\conv T$ cannot share two edges, since otherwise we would get a contradiction to \eqref{p in core}. We define $q, q_1, q_2, q_3, q_4$ such that $T = \{q, q_2, q_3 \},$ and $q_1,\ldots,q_4$ are consecutive vertices of $\conv A.$ Since $\card{A} \ge 6,$ for some $i=\{2,3\}$ the triangle $\conv \{q, q_{i+1},q_{i-1} \}$ does not share edges with $\conv A.$ Then $p \in \conv \{q,q_{i+1},q_{i-1} \}$, since otherwise one would get a contradiction to \eqref{p in core}. Hence for $T' := \{q, q_{i+1},q_{i-1} \}$, one has $p \in \conv (A \setminus T') \cap \conv T' \cap M,$ which shows that $r(M) \le h(M)+1.$ 
\end{proof}


\begin{thebibliography}{AKMM02}

\bibitem[AKMM02]{MR1921545}
N.~Alon, G.~Kalai, J.~Matou{\v{s}}ek, and R.~Meshulam, \emph{Transversal
  numbers for hypergraphs arising in geometry}, Adv. in Appl. Math. \textbf{29}
  (2002), no.~1, 79--101. \MR{2003g:52004}

\bibitem[ALW08]{CertificatesMixed2008}
K.~Andersen, Q.~Louveaux, and R.~Weismantel, \emph{Certificates of linear mixed
  integer infeasibility}, Operations Research Letters \textbf{36} (2008), 734
  -- 738.

\bibitem[BB03]{MR2007959}
K.~Bezdek and A.~Blokhuis, \emph{The {R}adon number of the three-dimensional
  integer lattice}, Discrete Comput. Geom. \textbf{30} (2003), no.~2, 181--184.
  \MR{2004g:52023}

\bibitem[Bel77]{Bell77}
D.~E. Bell, \emph{A theorem concerning the integer lattice}, Studies in Appl.
  Math. \textbf{56} (1977), no.~2, 187--188. \MR{57 \#2590}

\bibitem[BM03]{bar-mat-03}
I.~B{\'a}r{\'a}ny and J.~Matou{\v{s}}ek, \emph{A fractional {H}elly theorem for
  convex lattice sets}, Adv. Math. \textbf{174} (2003), no.~2, 227--235.
  \MR{2003m:52006}

\bibitem[BW05]{BertWeiBook}
D.~Bertsimas and R.~Weismantel, \emph{Optimization over {I}ntegers, {D}ynamic
  {I}deas}, Belmont, MA, 2005.

\bibitem[Cla95]{Clarkson95}
K.~L. Clarkson, \emph{Las {V}egas algorithms for linear and integer programming
  when the dimension is small}, J. Assoc. Comput. Mach. \textbf{42} (1995),
  no.~2, 488--499. \MR{97f:90041}

\bibitem[DGK63]{MR0157289}
L.~Danzer, B.~Gr{\"u}nbaum, and V.~Klee, \emph{Helly's theorem and its
  relatives}, Proc. {S}ympos. {P}ure {M}ath., {V}ol. {VII}, Amer. Math. Soc.,
  Providence, R.I., 1963, pp.~101--180. \MR{28 \#524}

\bibitem[Doi73]{MR0387090}
J.-P. Doignon, \emph{Convexity in crystallographical lattices}, J. Geometry
  \textbf{3} (1973), 71--85. \MR{52 \#7937}

\bibitem[Eck93]{MR1242986}
J.~Eckhoff, \emph{Helly, {R}adon, and {C}arath\'eodory type theorems}, Handbook
  of convex geometry, {V}ol.\ {A}, {B}, North-Holland, Amsterdam, 1993,
  pp.~389--448. \MR{94k:52010}

\bibitem[GPW93]{MR1228043}
J.~E. Goodman, R.~Pollack, and R.~Wenger, \emph{Geometric transversal theory},
  New trends in discrete and computational geometry, Algorithms Combin.,
  vol.~10, Springer, Berlin, 1993, pp.~163--198. \MR{95c:52010}

\bibitem[JW81]{MR634767}
R.~E. Jamison-Waldner, \emph{Partition numbers for trees and ordered sets},
  Pacific J. Math. \textbf{96} (1981), no.~1, 115--140. \MR{83a:52003}

\bibitem[KL79]{MR532152}
M.~Katchalski and A.~Liu, \emph{A problem of geometry in {${\bf R}\sp{n}$}},
  Proc. Amer. Math. Soc. \textbf{75} (1979), no.~2, 284--288. \MR{80h:52010}

\bibitem[Len83]{lenstra83}
H.~W. Lenstra, Jr., \emph{Integer programming with a fixed number of
  variables}, Math. Oper. Res. \textbf{8} (1983), no.~4, 538--548.
  \MR{86f:90106}

\bibitem[Onn91]{onn91}
S.~Onn, \emph{On the geometry and computational complexity of {R}adon
  partitions in the integer lattice}, SIAM J. Discrete Math. \textbf{4} (1991),
  no.~3, 436--446. \MR{92d:52041}

\bibitem[Roc97]{MR1451876}
R.~T. Rockafellar, \emph{Convex {A}nalysis}, Princeton University Press,
  Princeton, NJ, 1997. \MR{97m:49001}

\bibitem[Sca77]{MR0452678}
H.~E. Scarf, \emph{An observation on the structure of production sets with
  indivisibilities}, Proc. Nat. Acad. Sci. U.S.A. \textbf{74} (1977), no.~9,
  3637--3641. \MR{56 \#10957}

\bibitem[Sch86]{schrijver}
A.~Schrijver, \emph{Theory of {L}inear and {I}nteger {P}rogramming}, John Wiley
  \& Sons Ltd., Chichester, 1986. \MR{88m:90090}

\bibitem[Sch93]{schneider}
R.~Schneider, \emph{Convex {B}odies: the {B}runn-{M}inkowski {T}heory},
  Cambridge University Press, Cambridge, 1993. \MR{94d:52007}

\bibitem[Sie77]{MR0514026}
G.~Sierksma, \emph{Relationships between {C}arath\'eodory, {H}elly, {R}adon and
  exchange numbers of convexity spaces}, Nieuw Arch. Wisk. (3) \textbf{25}
  (1977), no.~2, 115--132. \MR{58 \#24020}

\bibitem[Tod77]{Todd1977}
M.~J. Todd, \emph{The number of necessary constraints in an integer program:
  {A} new proof of {S}carf's theorem}, Technical Report 355, School of
  Operations Research and Industrial Ingeneering, Cornell University, Ithaca,
  N.Y., 1977.

\bibitem[vdV93]{MR1234493}
M.~L.~J. van~de Vel, \emph{Theory of {C}onvex {S}tructures}, North-Holland
  Publishing Co., Amsterdam, 1993. \MR{95a:52002}

\bibitem[Zie95]{ziegler95}
G.~M. Ziegler, \emph{Lectures on {P}olytopes}, Springer-Verlag, New York, 1995.
  \MR{96a:52011}

\end{thebibliography}
\providecommand{\bysame}{\leavevmode\hbox to3em{\hrulefill}\thinspace}
\providecommand{\MR}{\relax\ifhmode\unskip\space\fi MR }
\providecommand{\MRhref}[2]{%
  \href{http://www.ams.org/mathscinet-getitem?mr=#1}{#2}
}
\providecommand{\href}[2]{#2}

\end{document}